\documentclass[12pt,oneside]{article}

\usepackage{amsthm,amsfonts,mathrsfs,  amsbsy, amssymb,amsmath,graphicx,euscript}

\usepackage{geometry}
\begin{document}

\begin{center} {\bf Solution-giving formula to  Cauchy problem\\ for multidimensional parabolic equation with variable coefficients

 Ivan\,D.\,Remizov}

National Research University Higher School of Economics

25/12 Bol. Pecherskaya Ulitsa, Room 224, Nizhny Novgorod, 603155, Russia

ivremizov@yandex.ru
\end{center}

We present a general method of solving the Cauchy problem for multidimensional parabolic (diffusion type) equation with variable coefficients which depend on spatial variable but do not change over time. We assume the existence of the $C_0$-semigroup (this is a standard assumption in the evolution equations theory, which guarantees the existence of the solution) and then find the representation (based on the family of translation operators) of the solution in terms of coefficients of the equation and initial condition. It is proved that if the coefficients of the equation are bounded, infinitely smooth and satisfy some other conditions then there exists a solution-giving $C_0$-semigroup of contraction operators.  We also represent the solution as a Feynman formula (i.e. as a limit of a multiple integral with multiplicity tending to infinity) with generalized functions appearing in the integral kernel.\\



\begin{center}
\textbf{Problem Setting and Approach Proposed}
\end{center}

At the present time we know of a relatively small number of situations where a short formula expresses the solution of a partial differential equation with variable coefficients in terms of these coefficients and initial/boundary conditions. In the present paper, we provide such formula for a diffusion-type equation (1), see below. Let us first describe the Cauchy problem and then provide the necessary background and references.

Consider integer dimension $d\geq 1$,  $x=(x_1,\dots,x_d)\in\mathbb{R}^d$, $t\geq0$, $u\colon [0,+\infty)\times \mathbb{R}^d\to \mathbb{R}$ and set the Cauchy problem for a second-order parabolic partial differential equation
\begin{equation}
\left\{ \begin{array}{ll}
   u'_{t}(t,x)=\sum_{j=1}^d(a_j(x))^2u_{x_jx_j}''(t,x)+\left<b(x), \nabla u(t,x)\right> +c(x)u(t,x)=Hu(t,x),  \\
   u(0,x)=u_0(x).
\end{array} \right.
\end{equation}

The coefficients of (1) are: an $\mathbb{R}^1$-valued function $c$, an $\mathbb{R}^d$-valued function $b$, and $\mathbb{R}^{1}$-valued functions $a_j$ for each $j=1,\dots,d$. We also represent function $b$ as a vector of $d$ $\mathbb{R}^1$-valued functions $b_j$ as follows: $b(x)=(b_1(x),\dots,b_d(x))$. We assume that all coefficients are bounded and uniformly continuous.   Parenthesis $\left<\cdot,\cdot\right>$ are used for the scalar product in $\mathbb{R}^d$. We use symbol $\nabla u(t,x)$ for the gradient vector with respect to $x$ as follows: $\nabla u=(\partial u/\partial x_1,\dots, \partial u/\partial x_d)$. For each $j=1,\dots,d$ we write the multiplier of higher derivative $\partial^2/\partial x_j^2$ in the form $(a_j(x))^2$ for two reasons: to show that the multiplier is non-negative, and to shorten some of the formulas. 

This paper is dedicated to deriving a formula that gives the solution of (1) in terms of $a_j$, $b_j$, $c$, $u_0$ assuming that the operator $H$ is an infinitesimal generator of the $C_0$-semigroup $\left(e^{tH}\right)_{t\geq 0}$. This assumption is standard in studies of evolution equations, which is the class of equations that the considered equation belongs to. According to  general theory of $C_0$-semigroups \cite{EN1} this assumption implies that the solution of the Cauchy problem (1) exists, is bounded and uniformly continuous with respect to $x$ for each $t$, depends on $u_0$ continuously, and can be represented in a form $u(t,x)=\left(e^{tH}u_0\right)(x)$. Similarly to what we have done in the case $d=1$ (see \cite{R-arch}) we apply the Chernoff theorem \cite{EN1, Chernoff} to a specially constructed family of operators $(S(t))_{t\geq 0}$, and express $e^{tH}$ in terms of $a_j$, $b_j$, $c$ reaching the proposed goal. We do not discuss the problem of finding the class of functions in which the solution is unique under certain assumptions on functions $a$, $b$, $c$, $u_0$, but keep in mind that for a heat equation there are known unbounded solutions. 

The formula that provides the solution of (1) is given in theorem 3.

\begin{center}
\textbf{State of the Art} 
\end{center}

The diffusion equation and heat equation have a long history of research since the beginning of 19-th century \cite{Histheat}, and there are many publications devoted to it.  The reader can see textbooks \cite{Ar1, Evans, HO}, recent papers \cite{BZu,OSY,PWX,jee1,jee2,jee3, Ar2, Ar3} and references therein, but still this list will be very incomplete. In the present paper, we are focused on the one  particular question: how we can express the solution of a parabolic (diffusion-type) equation with variable coefficients in terms of these coefficients? There are not so many publications answering it, but in the last decade the number of them increseases rapidly. This growth is achieved mainly via an approach based on the Chernoff theorem (theorem 1 below), and we also follow this way. Let us explain the idea of our method and mention some of the papers where the Chernoff theorem was used to obtain formulas (representations) for solutions of evolution equations.  

It is known (and will be discussed with more details in the ``Technique Employed'' section below), that if operator $H$ generates a $C_0$-semigroup $\left(e^{tH}\right)_{t\geq 0}$, then the solution of Cauchy problem $[u'_t(t,x)=Hu(t,x); u(0,x)=u_0(x)]$ can be represented in the form $u(t,x)=\left(e^{tH}u_0\right)(x)$. The Chernoff theorem \cite{EN1, Chernoff} allows us to reduce the problem of finding $e^{tH}$ to the problem of finding an appropriate operator-valued function $S(t)=I+tH+o(t)$, which is called the Chernoff function, and then use the Chernoff formula $e^{tH}=\lim_{n\to\infty}S(t/n)^n$. One advantage of that step is that we can define $S(t)$ by an explicit formula that depends on the coefficients of the operator $H$. Another advantage is that for each $t$ the operator $S(t)$ is a linear bounded operator, which  allows us to define analytic functions of argument $S(t)$ via power series (see examples \cite{R2,R3,R4,Sak,R-DAN2017}) to obtain a semigroup $\left(e^{-itH}\right)_{t\geq 0}$ that  solves Schr\"{o}dinger equation with Hamiltonian $H$.  This idea was introduced in \cite{R2} where we defined $R(t)=\exp\big(-i(S(t)-I)\big)$ and proved that $e^{-itH}=\lim_{n\to\infty}R(t/n)^n$. Members of O.G.Smolyanov's group employed Chernoff's theorem using integral operators as Chernoff functions to find solutions to parabolic equations in many cases during the last 15 years: see the pioneering papers  \cite{SWW, STT}, overview \cite{SmHist}, several examples \cite{SWW-pot, Butko1, R1, Dubravina, SmSh, Shamarov, bsg-2011} to see the diversity of applications, and recent papers \cite{R3, Buz, Sak, BGS-jmp, jee-Butko, OSS-randop, b-sd, R-DAN2017, OhShSh, KSS2016, Maz, Smolyanova2015, SKh2017, Burk2015, R-IDAQP-2018}. The solutions obtained were represented in the form of a \textit{Feynman formula}, i.e. as a limit of a multiple integral as the multiplicity goes to infinity. Indeed, if $S(t)$ is an integral operator for each $t>0$, then $S(t/n)^n$ is an $n-$tuple integral operator and equlity $e^{tH}u_0=\lim_{n\to\infty}S(t/n)^nu_0$ is a Feynman formula. See also \cite{R2} for applications of so-called \textit{quasi-Feynman formulas} which allow to solve Schr\"odinger equation in situations where it is difficult to obtain a Feynman formula for solution. See also papers \cite{M1, M2} and references therein, and Chapter 11 in \cite{JL}.

The specific feature of the presented research is that we use  translation operators instead of integral operators when constructing the Chernoff function $S(t)$, see \cite{BGS2010} for integral-based approach to similar equation and \cite{Plya1, Plya2} for equations with time-dependent coefficients; however, we allow all  $a_j(x)$ to be zero for some $x$, this case is not covered by methods of  \cite{BGS2010}. For this reason the solution of (1) is now represented via a new type of formulas that do not include integrals. However, one may interpret these formulas as Feynman formulas with Dirac delta-functions in the integral kernel, which is discussed later in remark 7 in the end of the paper.

\newpage
\begin{center}
	\textbf{Technique Employed} 
\end{center}

Let $\mathcal{F}$ be a Banach space. Let $\mathscr{L}(\mathcal{F})$ be a set of all bounded linear operators in $\mathcal{F}$. Suppose we have a mapping $V\colon [0,+\infty)\to \mathscr{L}(\mathcal{F}),$ i.e. $V(t)$ is a bounded linear operator $V(t)\colon \mathcal{F}\to \mathcal{F}$ for each $t\geq 0.$ The mapping $V$ is called \cite{EN1} a \textit{$C_0$-semigroup}, or \textit{a strongly continuous one-parameter semigroup} if it satisfies the following conditions: 
	
	1) $V(0)$ is the identity operator $I$, i.e. $\forall \varphi\in \mathcal{F}: V(0)\varphi=\varphi;$ 
	
	2) $V$ maps the addition of numbers in $[0,+\infty)$ into the composition of operators in $\mathscr{L}(\mathcal{F})$, i.e. $\forall t\geq 0,\forall s\geq 0: V(t+s)=V(t)\circ V(s),$ where for each $\varphi\in\mathcal{F}$ the notation $(A\circ B)(\varphi)=A(B(\varphi))=AB\varphi$ is used;
	
	3) $V$ is continuous with respect to the strong operator topology in $\mathscr{L}(\mathcal{F})$, i.e. $\forall \varphi\in \mathcal{F}$ function $t\longmapsto V(t)\varphi$ is continuous as a mapping $[0,+\infty)\to \mathcal{F}.$
	
	The definition of a \textit{$C_0$-group} is obtained by the substitution of $[0,+\infty)$ by $\mathbb{R}$ in the paragraph above.

It is known \cite{EN1} that if $(V(t))_{t\geq 0}$ is a $C_0$-semigroup in Banach space $\mathcal{F}$, then the set $$\left\{\varphi\in \mathcal{F}: \exists \lim_{t\to +0}\frac{V(t)\varphi-\varphi}{t}\right\}\stackrel{denote}{=}Dom(L)$$ is dense in $\mathcal{F}$. The operator $L$ defined on the domain $Dom(L)$ by the equality $$L\varphi=\lim_{t\to +0}\frac{V(t)\varphi-\varphi}{t}$$ is called \textit{an infinitesimal generator} (or just \textit{generator} to make it shorter) of the $C_0$-semigroup $(V(t))_{t\geq 0}$. 

One of the reasons for the study of $C_0$-semigroups is their connection with differential equations. If $Q$ is a set, then the function $u\colon [0,+\infty)\times Q\to \mathbb{C}$, $u\colon (t,x)\longmapsto u(t,x)$ of two variables $(t,x)$ can be considered as a function $u\colon t\longmapsto [x\longmapsto u(t,x)]$ of one variable $t$
with values in the space of functions of the variable $x$. If $u(t,\cdot)\in\mathcal{F}$ then one can define $Lu(t,x)=(Lu(t,\cdot))(x).$ If there exists a $C_0$-semigroup $(e^{tL})_{t\geq 0}$ then the Cauchy problem for a linear evolution equation 
$$
\left\{ \begin{array}{ll}
u'_t(t,x)=Lu(t,x) \ \mathrm{ for }\ t>0, x\in Q\\
u(0,x)=u_0(x)\ \mathrm{ for } \ x\in Q
\end{array} \right.
$$
has a unique (in sense of $\mathcal{F}$, where $u(t,\cdot)\in\mathcal{F}$ for every $t\geq 0$) solution $$u(t,x)=(e^{tL}u_0)(x)$$ depending on $u_0$ continuously. Compare also different meanings of the solution \cite{EN1}, including mild solution which solves the corresponding integral equation. Note that if there exists a strongly continuous group $(e^{tL})_{t\in\mathbb{R}}$ then in the Cauchy problem the equation $u'_t(t,x)=Lu(t,x)$ can be considered not only for $t>0$, but for $t\in\mathbb{R}$, and the solution is provided by the same formula $u(t,x)=(e^{tL}u_0)(x)$.

\textbf{Definition 1} (\textit{Introduced in \cite{R2}}). Let us say that $G$ is \textit{Chernoff-tangent} to $L$ iff the following conditions of Chernoff tangency (CT) hold: 

(CT0). Let $\mathcal{F}$ be a Banach space, and $\mathscr{L}(\mathcal{F})$ be a space of all linear bounded operators in $\mathcal{F}$. Suppose that we have an operator-valued function $G\colon [0, +\infty) \to \mathscr{L}(\mathcal{F})$, or, using other words, we have a family $(G(t))_{t\geq 0}$ of linear bounded operators in $\mathcal{F}$. Closed linear operator $L\colon Dom(L) \to \mathcal{F}$ is defined on the linear subspace $Dom(L)\subset\mathcal{F}$ which is dense in  $\mathcal{F}$.

(CT1) Function $t\longmapsto G(t)f\in\mathcal{F}$ is continuous for each $f\in\mathcal{F}$. 

(CT2) $G(0)=I$, i.e. $G(0)f=f$ for each $f\in\mathcal{F}$.

(CT3) There exists such a dense subspace $\mathcal{D}\subset \mathcal{F}$ that for each $f\in \mathcal{D}$ there exists a limit $$G'(0)f=\lim_{t\to 0}\frac{G(t)f-f}{t}.$$ 

(CT4) The closure of the operator  $(G'(0),\mathcal{D})$ is equal to $(L,Dom(L)).$ 

\textbf{Remark 1.} Let us consider one-dimensional example $\mathcal{F}=\mathscr{L}(\mathcal{F})=\mathbb{R}$. Then $g\colon [0,+\infty)\to\mathbb{R}$ is Chernoff-tangent to $l\in\mathbb{R}$ iff $g(t)=1+tl+o(t)$ as $t\to+0$.

\textbf{Theorem 1} (\textsc{P.\,R.~Chernoff (1968)}, see \cite{EN1, Chernoff}). Let $\mathcal{F}$ and $\mathscr{L}(\mathcal{F})$ be as above. Suppose that the operator $L\colon \mathcal{F}\supset Dom(L)\to \mathcal{F}$ is linear and closed, and function $G$ takes values in $\mathscr{L}(\mathcal{F})$. Suppose that these assumptions are fulfilled:

(E) There exists a $C_0$-semigroup $(e^{tL})_{t\geq 0}$ with the infenitesimal generator $(L,Dom(L))$.

(CT) $G$ is Chernoff-tangent to $(L,Dom(L)).$ 

(N) There exists such a number $\omega\in\mathbb{R}$, that $\|G(t)\|\leq e^{\omega t}$ for all $t\geq 0$.

Then for each $f\in \mathcal{F}$  we have $(G(t/n))^nf\to e^{tL}f$ as $n\to \infty$ with respect to norm in $\mathcal{F}$ uniformly with respect to $t\in[0,T]$ for each $T>0$, i.e.
$$\lim_{n\to\infty}\sup_{t\in[0,T]}\left\|e^{tL}f - (G(t/n))^nf \right\|=0.$$

\textbf{Remark 2.} In our one-dimensional example ($\mathcal{F}=\mathscr{L}(\mathcal{F})=\mathbb{R}$) the Chernoff theorem says that  $e^{tl}=\lim_{n\to\infty}g(t/n)^n=\lim_{n\to\infty}(1+tl/n+o(t/n))^n$, which is a simple fact of calculus.

\textbf{Definition 2.} Let $\mathcal{F}, \mathscr{L}(\mathcal{F}), L$ be as above. If $G$ is Chernoff-tangent to $L$ and the equation $\lim_{n\to\infty}\sup_{t\in[0,T]}\left\|e^{tL}f - (G(t/n))^nf \right\|=0$ holds, then $G$ is called \textit{a Chernoff function} for the operator $L$, and the $(G(t/n))^nf$ is called \textit{a Chernoff approximation expression} to $e^{tL}f$. 

\textbf{Remark 3.} If $L$ is a linear bounded operator in $\mathcal{F}$, then $e^{tL}=\sum_{k=0}^{+\infty}(tL)^k/k!$ where the series converges in the usual operator norm topology in $\mathscr{L}(\mathcal{F})$. When $L$ is not bounded (such as Laplacian and many other differential operators), expressing $(e^{tL})_{t\geq 0}$ in terms of $L$ is not an easy problem that is equivalent to the problem of finding (for each $u_0\in\mathcal{F}$) the $\mathcal{F}$-valued function $U$ that solves the Cauchy problem $U'(t)=LU(t); U(0)=u_0$. If one finds this solution, then $e^{tL}$ is  obtained for each $u_0\in\mathcal{F}$ and each $t\geq 0$ in the form $e^{tL}u_0=U(t)$.

\textbf{Remark 4.} In the definition of the Chernoff tangency the family $(G(t))_{t\geq 0}$ usually does not have a semigroup composition property, i.e. $G(t_1+t_2)\neq G(t_1)G(t_2)$,  while $(e^{tL})_{t\geq 0}$ has it: $e^{t_1L}e^{t_2L}=e^{(t_1+t_2)L}$. However, each $C_0$-semigroup $(e^{tL})_{t\geq 0}$ is Chernoff-tangent to its generator $L$ and appears to be it's Chernoff function. When coefficients of the operator $L$ are variable, usually there is no simple formula for $e^{tL}$ due to the remark 3. On the other hand, even in this case one can find rather simple formula to construct Chernoff function $G$ for the operator $L$, because there is no need to worry about the composition property, and then obtain $e^{tL}$ in the form $e^{tL}=\lim_{n\to\infty}G(t/n)^n$ via the Chernoff theorem. This is what we do in the present paper for $L=H$ and what people have done for different operators described in papers cited above.

\textbf{
\begin{center}
Chernoff Function for Operator $H$
\end{center}
}

\textbf{Remark 5.} Let us denote the set of all (real-valued and defined on $\mathbb{R}^d$) bounded continuous functions as $C_b(\mathbb{R}^d)$, the set of all bounded functions with bounded derivatives of all orders as $C^{\infty}_b(\mathbb{R}^d)$, and the set of all bounded, uniformly continuous functions as $UC_b(\mathbb{R}^d)$. Then $C^{\infty}_b(\mathbb{R}^d)\subset UC_b(\mathbb{R}^d)\subset C_b(\mathbb{R}^d)$, and with respect to the uniform (Chebyshev) norm $\|f\|=\sup_{x\in\mathbb{R}^d}|f(x)|$ the first inclusion is dense, and the last two spaces are Banach spaces. 

Indeed: the inclusions follow directly from the definitions; $C_b(\mathbb{R}^d)$ is a Banach space -- standard fact; $UC_b(\mathbb{R}^d)$ is a closed subset of $C_b(\mathbb{R}^d)$ -- follows from the fact that the uniform limit of a sequence of uniformly continuous bounded functions is a uniformly continuous bounded function, simple check; hence $UC_b(\mathbb{R}^d)$ is a Banach space being a closed subset of a Banach space $C_b(\mathbb{R}^d)$; uniform closure of $C^{\infty}_b(\mathbb{R}^d)$ is included in $UC_b(\mathbb{R}^d)$ -- because $C^{\infty}_b(\mathbb{R}^d)\subset UC_b(\mathbb{R}^d)$ and $UC_b(\mathbb{R}^d)$ is closed; finally, let us prove that uniform closure of $C^{\infty}_b(\mathbb{R}^d)$ includes $UC_b(\mathbb{R}^d)$.

\textbf{Lemma 1.} For each function $f\in UC_b(\mathbb{R}^d)$ and each $\varepsilon>0$ there exists a function $\phi\in C^{\infty}_b(\mathbb{R}^d)$ such that $\|f-\phi\|\leq \varepsilon$.
 
\textbf{Proof.} Suppose that $f\in UC_b(\mathbb{R}^d)$ and $\varepsilon>0$ are given. We will find such a $\psi\in C^{\infty}_c(\mathbb{R}^d)$  that for $\phi=f*\psi$ relation  $\|f-\phi\|\leq\varepsilon$ holds, where symbol  $C^{\infty}_c(\mathbb{R}^d)$ denotes the space of all real-valued infinitely differentiable functions on $\mathbb{R}^d$ with compact support --- this means that each of them vanishes outside some (depending on the function) ball in $\mathbb{R}^d$, is bounded and all the derivatives are also bounded.

a) It is known (see \cite{Chist}, Proposition 1.3 and notation above it) that for each $f\in UC_b(\mathbb{R}^d)\subset L_1^{loc}(\mathbb{R}^d)$ and $\psi\in C^{\infty}_c(\mathbb{R}^d)$ function $f*\psi$ is continuous, has derivatives of all orders and for each multi-index $\alpha$ we have $\partial^\alpha (f*\psi)=f*(\partial^\alpha\psi).$ Note that $\psi\in C^{\infty}_c(\mathbb{R}^d)$ implies $\partial^\alpha\psi\in C^{\infty}_c(\mathbb{R}^d)$.

b) Let us prove that the convolution of $f$ with any $\eta\in C^{\infty}_c(\mathbb{R}^d)$ is a bounded function; due to a) this will give us that function $f*\psi$ is bounded and all its derivatives are also bounded. We know that $\eta$ iz zero outside some ball   $B\subset\mathbb{R}^d$. Then for all $x\in\mathbb{R}^d$ we have  $|(f*\eta)(x)|=\left|\int_{\mathbb{R}^d}f(x-y)\eta(y)dy\right|=\left|\int_{B}f(x-y)\eta(y)dy\right|\leq \sup_{z\in\mathbb{R}^d}|f(z)|\cdot \sup_{y\in B}|\eta(y)|\cdot vol(B)<\infty$. So $(f*\psi)\in C^{\infty}_b(\mathbb{R}^d)$ for all $\psi\in C^{\infty}_c(\mathbb{R}^d)$.

c) For $x\in\mathbb{R}^d$ employ notation $|x|=\sqrt{x_1^2+\dots+x_d^2}$ and for each $\delta>0$ define $$\varphi_\delta(x)=\exp\left(\frac{-1}{1-|x/\delta|^2}\right)\cdot \left(\int_{|y|<\delta}\exp\left(\frac{-1}{1-|y/\delta|^2}\right)dy\right)^{-1} \textrm{ for }|x|<\delta,\quad \varphi_\delta(x)=0 \textrm{ for }|x|\geq\delta.$$
Simple check shows that $\varphi_\delta\in C^{\infty}_c(\mathbb{R}^d)$, $\int_{\mathbb{R}^d}\varphi_\delta(x)dx=1$, $0\leq \varphi(x)\leq \max_{x\in\mathbb{R}^d}\varphi_\delta(x)\stackrel{denote}{=}M(\delta)<\infty$ where $\lim_{\delta\to0}M(\delta)=+\infty$.

d) As $f$ is uniformly continuous, there exists such $\delta>0$ that $|y|<\delta$ implies $|f(x-y)-f(x)|<\varepsilon$ for all $x\in\mathbb{R}^d$.

e) Define $\psi=\varphi_\delta$ and $\phi(x)=(f*\psi)(x)=\int_{\mathbb{R}^d}f(x-y)\psi(y)dy$. 

f) The estimate $|\phi(x)- f(x)\cdot 1|=\left|\int_{\mathbb{R}^d}f(x-y)\psi(y)dy-f(x)\int_{\mathbb{R}^d}\psi(y)dy\right| \leq \int_{\mathbb{R}^d}|f(x-y)-f(x)|\psi(y)dy=\int_{|y|<\delta}|f(x-y)-f(x)|\psi(y)dy< \int_{|y|<\delta}\varepsilon\psi(y)dy=\varepsilon$ holds for all $x\in\mathbb{R}^d$ so $\|\phi-f\|=\sup_{x\in\mathbb{R}^d}|\phi(x)- f(x)|\leq \varepsilon$. $\Box$

\textbf{Theorem 2.} Let $e_j\in\mathbb{R}^d$ be a constant $d$-dimensional vector with $1$ at position $j$ and $0$ at other $d-1$ positions. For each $x\in\mathbb{R}^d$, $t\geq 0$, $f\in C_b(\mathbb{R}^d)$ and $\varphi\in C^{\infty}_b(\mathbb{R}^d)$ set
\begin{multline}(S(t)f)(x)=\frac{1}{4d}\sum_{j=1}^d\bigg(f\left(x+2\sqrt{d}a_j(x)\sqrt{t}e_j\right)+\\
+f\left(x-2\sqrt{d}a_j(x)\sqrt{t}e_j\right)\bigg)+\frac{1}{2}f(x+2tb(x))+tc(x)f(x),\end{multline}

$$(H\varphi)(x)=
\sum_{j=1}^d(a_j(x))^2\varphi_{x_jx_j}''(x)+\left<b(x), \nabla \varphi(x)\right> +c(x)\varphi(x).\eqno(3)$$
Then, with respect to the norm $\|g\|=\sup_{x\in\mathbb{R}^d}|g(x)|$, the following holds:

I) for each $t\geq 0$ and $f\in C_b(\mathbb{R})$ we have $\|S(t)f\|\leq \big(1+\|c\|t\big) \|f\|$.

II) for each $\varphi\in C^{\infty}_b(\mathbb{R}^d)$ we have $\lim_{t\to+0}\|S(t)\varphi-\varphi-tH\varphi\|/t=0$.

III) if $t_n\to t_0$, $t_n\geq 0$ and $f\in UC_b(\mathbb{R}^d)$, then $\lim\limits_{t\to t_0}\|S(t_n)f- S(t_0)f\|=0$ for each $t_0\geq 0$.
 
IV) if $a_j,b_j,c,f\in UC_b(\mathbb{R}^d)$, then $S(t)f\in UC_b(\mathbb{R}^d)$ for each $t\geq 0$.

\textbf{Proof.} I) Let us write $\sup$ instead of $\sup_{x\in\mathbb{R}^d}$ in the proof of this item to make it shorter.  

$\|S(t)f\|=\sup\Big|\frac{1}{4d}\sum_{j=1}^d f\left(x+2\sqrt{d}a_j(x)\sqrt{t}e_j\right)+\frac{1}{4d}\sum_{j=1}^df\left(x-2\sqrt{d}a_j(x)\sqrt{t}e_j\right)+\frac{1}{2}f(x+2tb(x))+tc(x)f(x)\Big| \leq \frac{1}{4d}\sum_{j=1}^d\sup\left|f\left(x+2\sqrt{d}a_j(x)\sqrt{t}e_j\right)\right|$

$+\frac{1}{4d}\sum_{j=1}^d\sup\left|f\left(x-2\sqrt{d}a_j(x)\sqrt{t}e_j\right)\right| +\frac{1}{2}\sup|f(x+2b(x)t)|+t\sup|c(x)|\sup|f(x)|\leq$

$\frac{1}{4d}d\|f\|+\frac{1}{4d}d\|f\|+\frac{1}{2}\|f\|+t\sup|c(x)|\|f\|=\big(1+\|c\|t\big) \|f\|.$

II) Now remember that function $\varphi$ in this item is bounded with all derivatives. Let us fix arbitrary $x\in\mathbb{R}$ in (2) and consider $(S(t)\varphi)(x)$ as a smooth function of $\sqrt{t}$. Then we use Taylor's expansion in powers of $\sqrt{t}$ in first two sums and in powers of $t$ in the third summand to show that $(S(t)\varphi)(x)=\varphi(x)+t(H\varphi)(x)+t\sqrt{t}R(t,x)$. As the derivatives of $\varphi$ are bounded, and functions $a_j,b_j,c$ are also bounded, one can represent the remainders in Lagrange's form and see that $\sup_{t\in [0,t_0]}\sup_{x\in\mathbb{R}^d}|R(t,x)|<\infty$ for each fixed $t_0>0$. We skip the exact formula for $R(t,x)$ to make the proof shorter. 
Indeed, for each $j$ we have  $$\varphi\left(x \pm  2\sqrt{d}a_j(x)\sqrt{t}e_j\right)=\varphi(x)\pm 2\sqrt{d}a_j(x)\sqrt{t}\left<\nabla \varphi(x),e_j\right> +\frac{1}{2}(2\sqrt{d}a_j(x)\sqrt{t})^2\left<\varphi''(x)e_j,e_j\right>+o(t)$$
so
$$\varphi\left(x +  2\sqrt{d}a_j(x)\sqrt{t}e_j\right)  + \varphi\left(x -  2\sqrt{d}a_j(x)\sqrt{t}e_j\right)=2 \varphi(x)+4d(a_j(x))^2\left<\varphi''(x)e_j,e_j\right>t +o(t),$$
and after summation we obtain
$$
\frac{1}{4d}\sum_{j=1}^d\left( \varphi\left(x +  2\sqrt{d}a_j(x)\sqrt{t}e_j\right)  + \varphi\left(x -  2\sqrt{d}a_j(x)\sqrt{t}e_j\right) \right)=\frac{1}{2}\varphi(x) +t\sum_{j=1}^d(a_j(x))^2\varphi_{x_jx_j}''(x) +o(t).
$$
For the third summand we have 
$$\frac{1}{2}\varphi(x+2tb(x))=\frac{1}{2}\varphi(x) +t\left<\nabla \varphi(x),b(x)\right>+o(t).$$
The term $tc(x)\varphi(x)$ is already in the form we need. Summing up we have
$$(S(t)\varphi)(x)=
\frac{1}{4d}\sum_{j=1}^d\left( \varphi\left(x +  2\sqrt{d}a_j(x)\sqrt{t}e_j\right)  + \varphi\left(x -  2\sqrt{d}a_j(x)\sqrt{t}e_j\right) \right) + \frac{1}{2}\varphi(x+2tb(x)) +tc(x)\varphi(x)
$$
$$
=\left(\frac{1}{2}\varphi(x) +t\sum_{j=1}^d(a_j(x))^2\varphi_{x_jx_j}''(x) +o(t) \right) + \left(\frac{1}{2}\varphi(x) +t\left<\nabla \varphi(x),b(x)\right>+o(t) \right) + tc(x)\varphi(x)
$$
$$
=\varphi(x) +t\sum_{j=1}^d(a_j(x))^2\varphi_{x_jx_j}''(x) + t\left<\nabla \varphi(x),b(x)\right> +tc(x)\varphi(x) +o(t)=\varphi(x)+t(H\varphi)(x)+o(t).
$$

III) In this item we assume that $t_n\to t_0$ and prove that $(S(t_n)f)(x)\to (S(t_0)f)(x)$ uniformly with respect to $x\in\mathbb{R}^d$ as $n\to\infty$. Indeed, functions $a_j$ are bounded, so for each $j$ we have $x +  2\sqrt{d}a_j(x)\sqrt{t_n}e_j\to x +  2\sqrt{d}a_j(x)\sqrt{t_0}e_j$ uniformly with respect to $x$. Function $f$ is uniformly continuous, so $f(x +  2\sqrt{d}a_j(x)\sqrt{t_n}e_j)\to f(x +  2\sqrt{d}a_j(x)\sqrt{t_0}e_j)$ uniformly with respect to $x$. In the same manner we use the fact that functions $b_j, c$ are bounded, and then sum all the limit conditions obtained.

IV) If $t\geq 0$ is fixed, then $\big[z\longmapsto f(z +  2\sqrt{d}a_j(z)\sqrt{t_n}e_j)\big]\in UC_b(\mathbb{R}^d)$ because  $a_j,f\in UC_b(\mathbb{R}^d)$. All the summands of $S(t)f$ are processed in this manner. $\Box$ 

\newpage
\textbf{
	\begin{center}
		Main Result
	\end{center}
}

\textbf{Theorem 3.} \textit{(On the representation of the solution of the Cauchy problem in case when the closure of $H$ generates a $C_0$-semigroup.)} Suppose that functions $a_j$, $b_j$, $c$ belong to the space $UC_b(\mathbb{R}^d)$ endowed with the norm $\|f\|=\sup_{x\in\mathbb{R}^d}|f(x)|$. Suppose that operator $H$ is defined by equation (3) on the domain $C^\infty_b(\mathbb{R}^d)\subset UC_b(\mathbb{R}^d)$, and the closure of this operator: a)~exists;
b)~is~an~infinitesimal generator of a $C_0$-semigroup $(e^{tH})_{t\geq0}$ in $UC_b(\mathbb{R}^d)$.

Then for each $u_0\in UC_b(\mathbb{R}^d)$ there exists a bounded (and uniformly continuous with respect to $x\in\mathbb{R}^d$ for each $t\geq 0$) solution $u$ of the Cauchy problem~(1), it depends on $u_0$ continuously and uniformly with respect to $x\in\mathbb{R}^d$ for each $t\geq 0$. For each $x\in\mathbb{R}^d$ and $t\geq0$ this solution is given by the formula $$u(t,x)=\left(e^{tH}u_0\right)(x)=\lim_{n\to\infty}\Big(\Big(S(t/n)\Big)^nu_0\Big)(x),\eqno(4)$$
where $S(t/n)$ is obtained by substitution of $t$ by $t/n$ in the equation (2), and the $n$-th power in the expression $(S(t/n))^n$ means the composition of $n$ copies of linear bounded translation operator $S(t/n)$. The limit (4) for each fixed $t>0$ is taken in the space $UC_b(\mathbb{R}^d)$ and appears to be uniform with respect to $t\in[0,t_0]$ for each $t_0>0$.

\textbf{Proof.} Let us check the conditions of the Chernoff theorem and thus show that $S$ is a Chernoff function for $L$. In theorem~1 and definition 1 we set $\mathcal{F}=UC_b(\mathbb{R})$, $G(t)=S(t)$, $L=H$, $\mathcal{D}=C^\infty_b(\mathbb{R})$, $\omega=\|c\|$. Condition $(E)$ is an assumption from theorem 3, condition $(N)$ is provided by item I) of  theorem 2: $\|S(t)\|\leq 1+\|c\|t\leq e^{\|c\|t}$. Let us check the Chernoff tangency: (CT0) follows from remark 5, lemma 1 and assumptions of theorem 3; (CT1) follows from items IV) and III) of theorem 2; (CT2) follows directly from formula (2); (CT3) follows from item II) of theorem 2; (CT4) is an assumption of theorem 3. Therefore the statement of theorem 3 is true thanks to the statement of the Chernoff theorem and standard facts of the $C_0$-semigroup theory \cite{EN1}. $\Box$ 

\textbf{Remark 6.} Formula (4) proven in theorem~3 contains $\lim_{n\to\infty}$. After the limit is taken we obtain the exact solution to Cauchy problem (1). For each fixed $n$ the expression under the limit sign is an approximation of the solution. With growth of $n$ such approximations converge to the exact solution uniformly with respect to $x\in\mathbb{R}^d$ and $t\in[0,t_0]$ for each fixed $t_0>0$.

\textbf{Lemma 2.} \textit{(On conditions that guarantee that the closure of $H$ generates a $C_0$-semigroup.)} Suppose that in (3) coefficients $a_j, b_j, c$ belong to space $C_b^\infty(\mathbb{R}^d)$,  $c(x)\leq 0$ and there exists a constant $\varkappa>0$ such that for each $\xi=(\xi_1,\dots,\xi_n)\in\mathbb{R}^d$ and all $x\in\mathbb{R}^d$
the ellipticity condition is fulfilled:  $\sum_{j=1}^da_{j}(x)^2\xi^2_j\geq\varkappa\sum_{j=1}^d\xi^2_j.$ Then:

I) For all $\lambda>0$ we have $C_b^\infty(\mathbb{R}^d)\subset (H-\lambda I)(C_b^\infty(\mathbb{R}^d))$.

II) Operator $H$ defined by (3) on domain $C_b^\infty(\mathbb{R}^d)$ is dissipative.

III) Closure of $(H, C_b^\infty(\mathbb{R}^d))$ is a dissipative operator that  generates a contraction $C_0$-semigroup in $UC_b(\mathbb{R}^d)$.

IV) Cauchy problem (1) has the solution $u(t,x)$ which is given by $u(t,x)=\left(e^{tH}u_0\right)(x)$ and for all $t\geq 0$ satisfies $\sup_{x\in\mathbb{R}^d}|u(t,x)|\leq \sup_{x\in\mathbb{R}^d}|u_0(x)|$.

\textbf{Proof}. Items I) and II) follow respectively from items (1) and (2) of Lemma 2.5 in \cite{R-IDAQP-2018}. Each dissipative operator is closable and its closure is again a dissipative operator -- item (iv) of Proposition 3.14 in \cite{EN1}. Let us note that due to Lemma 1 item I) implies that the range of operator $\lambda I-H$ is dense in $UC_b(\mathbb{R}^d)$ for all $\lambda>0$. Item III) now follows from this fact, item II and theorem 3.15 in \cite{EN1} which is also known as Lumer-Philipps theorem. Item IV) follows from item III), Proposition 6.2 in \cite{EN1} and the definition of a contraction operator which in this case says that $\left\|e^{tH}\right\|\leq 1$ for all $t\geq 0$. Let us mention also that we have just proved that under assumptions of lemma 2 the space $C_b^\infty(\mathbb{R}^d)$ is a core for operator $H$. $\Box$

\textbf{Remark 7.} Let us note that assumptions of lemma 2 may be relaxed in several ways, see remark 3.1 in \cite{BGS2010}. Considering spaces of functions vanishing at infinity and combining ideas of \cite{BGS2010} with the technique of the present paper one can try to yield representation of solutions to initial-value problems for PDEs with not so smooth coefficients, and of solutions to initial-boundary problems in bounded domains in $\mathbb{R}^d$ with smooth (but not infinitely smooth) boundary.

\textbf{
\begin{center}
Feynman Formulas with Generalized Functions
\end{center}
}
\textbf{Remark 8.} Formula (2) can be rewritten in terms of generalized functions (=distributions) in $\mathbb{R}^d$, based on the fact that for Dirac's $\delta$-function the following equation holds by definition of the integral in the right-hand side:  $f(w)=\int_{\mathbb{R}^d}\delta(y-w)f(y)dy$. We fix $x\in\mathbb{R}^d$ and  start from
$$f\left(x+2\sqrt{d}a_j(x)\sqrt{t}e_j\right)=\int_{\mathbb{R}^d}\delta\left(y-x-2\sqrt{d}a_j(x)\sqrt{t}e_j\right)f(y)dy,$$
then, after the substitution $y=x+z,  z=y-x, dy=dz$ in the integral above, we get
$$f\left(x+2\sqrt{d}a_j(x)\sqrt{t}e_j\right)=\int_{\mathbb{R}^d}\delta\left(z-2\sqrt{d}a_j(x)\sqrt{t}e_j\right)f(x+z)dy.$$
In the same way we have
$$f\left(x-2\sqrt{d}a_j(x)\sqrt{t}e_j\right)=\int_{\mathbb{R}^d}\delta\left(z+2\sqrt{d}a_j(x)\sqrt{t}e_j\right)f(x+z)dy,$$
$$f(x+2tb(x))=\int_{\mathbb{R}^d}\delta\left(z-2tb(x)\right)f(x+z)dy,\quad  c(x)f(x)=\int_{\mathbb{R}^d}c(x)\delta(z)f(x+z)dy.$$
Finally
$$(S(t)f)(x)=\int_{\mathbb{R}^d} \Phi(z,x,t)f(x+z)dz,\eqno(5)$$
where
$$
\Phi(z,x,t)=\frac{1}{4d}\sum_{j=1}^d\left(\delta\left(z-2\sqrt{d}a_j(x)\sqrt{t}e_j\right) + \delta\left(z+2\sqrt{d}a_j(x)\sqrt{t}e_j\right)\right)+\frac{1}{2}\delta\left(z-2tb(x)\right)+tc(x)\delta(z).
$$
Then (4) can be rewritten as a Feynman formula, i.e. as a representation of the function $u$ in a form of a limit of multiple integral where multiplicity tends to infinity:
$$u(t,x)=\lim_{n\to\infty}(S(t/n)^nu_0)(x)
=
\lim_{n\to\infty}\int_{\mathbb{R}^d}\Phi(z_1,x, t/n)\int_{\mathbb{R}^d}\Phi(z_2,x+z_1, t/n)\int_{\mathbb{R}^d}\Phi(z_3,x+z_1+z_2, t/n)\dots$$
$$\dots\int_{\mathbb{R}^d}\Phi(z_n,x+z_1+\dots+z_{n-1}, t/n)u_0(x+z_1+\dots+z_n)dz_n\dots dz_1.\eqno(6)$$

\textbf{Remark 9.} Right-hand sides of equations (5) and (6) are formal expressions, but left-hand sides are well-defined thanks to theorem 2 and theorem 3. It is important to note that multiple integrals in (6) are not true multiple integrals that can be rewritten as one integral over $\mathbb{R}^{nd}$. Instead (6) needs to be understood as several integrals that are taken one after another: first $dz_n$, then $dz_{n-1}$ and finishing with $dz_1$. Indeed, if we write an expression
$$\lim_{n\to\infty}
\underbrace{\int_{\mathbb{R}^d}\dots  
	\int_{\mathbb{R}^d}}_n
\Phi(z_1,x,t/n)\dots\Phi(z_n,x+z_1+\dots+z_{n-1}, t/n)u_0(x+z_1+\dots+z_n)dz_n\dots dz_1\eqno(7)$$
then under the integral sign we will have (tensor in this case) product of delta-functions which always needs a delicate handling. However, it seems interesting to study such expressions as (7) using methods of generalized functions theory \cite{Chist}. It is also interesting to find the physical meaning of generalized function $\Phi$ which in regular case is a transitional density of the diffusion process. Another way to write the same is provided by dual paring notation, but in this notation analogy with Feynman formulas is not very well seen.

\textbf{
\begin{center}
Acknowledgements
\end{center}
}
The author is thankful to his scientific advisor Prof. O.G.Smolyanov for encouragement and consultations, to Prof.  D.V.Turaev for delightful discussions and to Prof. S.Mazzucchi for helpful references. It is also a great plesure to thank the referee of J. Math. Phys. who opened several important questions (which were answered in lemma 1, lemma 2, remark 7, remark 9) and Prof. V.V.Chistyakov who provided a great reference paper \cite{Chist} which helped in proving lemma 1. This research was carried out within the HSE University Basic Research Program in 2019.


\end{document}